\documentclass{article}
\usepackage{amssymb}
\usepackage{amsmath}

\newtheorem{theorem}{Theorem}

\newtheorem{corollary}[theorem]{Corollary}

\newtheorem{lemma}[theorem]{Lemma}

\newtheorem{proposition}[theorem]{Proposition}

\begin{document}

\title{The fractional Poisson measure in infinite dimensions}
\author{M.~J.~Oliveira\thanks{%
Centro de Matem\'{a}tica e Aplica\c{c}\~{o}es Fundamentais, Av.~Prof.~Gama
Pinto 2, 1649-003 Lisboa, Portugal} \thanks{%
Universidade Aberta, Lisboa, Portugal, oliveira@cii.fc.ul.pt}, H.~Ouerdiane%
\thanks{%
D{\'{e}}partement de Math{\'{e}}matiques, Facult{\'{e}} des Sciences de
Tunis, 1060 Tunis, Tunisia, habib.ouerdiane@fst.rnu.tn}, J.~L.~Silva\thanks{%
CCM, University of Madeira, 9000-390 Funchal, Portugal, luis@uma.pt} and
R.~Vilela Mendes\footnotemark[1] \thanks{%
Instituto de Plasmas e Fus\~{a}o Nuclear, IST, Lisboa, Portugal,
rvilela.mendes@gmail.com, corresponding author}}
\maketitle

\begin{abstract}
The Mittag-Leffler function $E_{\alpha }$ being a natural generalization of
the exponential function, an infinite-dimensional version of the fractional
Poisson measure would have a characteristic functional 
\begin{equation*}
C_{\alpha }\left( \varphi \right) :=E_{\alpha }\left( \int (e^{i\varphi
\left( x\right) }-1)\,d\mu \left( x\right) \right)
\end{equation*}%
which we prove to fulfill all requirements of the Bochner-Minlos theorem.

The identity of the support of this new measure with the support of the
infinite-dimensional Poisson measure ($\alpha =1$) allows the development of
a fractional infinite-dimensional analysis modeled on Poisson analysis
through the combinatorial harmonic analysis on configuration spaces. This
setting provides, in particular, explicit formulas for annihilation,
creation, and second quantization operators. In spite of the identity of the
supports, the fractional Poisson measure displays some noticeable
differences in relation to the Poisson measure, which may be physically
quite significant.

Keywords: Poisson measure, Fractional, Infinite-dimensional analysis

MSC 28C20, 60G55
\end{abstract}

\section{Introduction}

The Poisson measure $\pi$ in $\mathbb{R}$\ (or $\mathbb{N}$) is 
\begin{equation*}
\pi \left( A\right) =e^{-\sigma }\sum_{n\in A}\frac{\sigma ^{n}}{n!}
\label{1.1}
\end{equation*}
the parameter $\sigma $ being called the \textit{intensity\/}. The Laplace
transform of $\pi$ is 
\begin{equation*}
l_{\pi }\left( \lambda \right) =\mathbb{E}\left( e^{\lambda \cdot }\right)
=e^{-\sigma }\sum_{n=0}^{\infty }\frac{\sigma ^{n}}{n!}e^{\lambda
n}=e^{\sigma \left( e^{\lambda }-1\right) }  \label{1.2}
\end{equation*}
For $n$-tuples of independent Poisson variables one would have 
\begin{equation*}
l_{\pi }\left( \boldsymbol{\lambda }\right) =e^{\sum \sigma _{k}\left(
e^{\lambda _{k}}-1\right) }  \label{1.3}
\end{equation*}
Continuing $\lambda _{k}$ to imaginary arguments $\lambda _{k}=if_{k}$, the
characteristic function is 
\begin{equation}
C_{\pi }\left( \lambda \right) =e^{\sum \sigma _{k}\left(
e^{if_{k}}-1\right) }  \label{1.4}
\end{equation}

Looked at as a renewal process, $P\left( X=n\right) =e^{-\sigma }\frac{%
\sigma ^{n}}{n!}$ would be the probability of $n$ events occurring in the
time interval $\sigma $. The survival probability, that is, the probability
of no event is 
\begin{equation*}
\Psi \left( \sigma \right) =e^{-\sigma }
\end{equation*}%
which satisfies the equation 
\begin{equation}
\frac{d}{d\sigma }\Psi \left( \sigma \right) =-\Psi \left( \sigma \right)
\label{1.6}
\end{equation}

Replacing in (\ref{1.6}) the derivative $\frac{d}{d\sigma }$ by the (Caputo)
fractional derivative 
\begin{equation*}
D^{\alpha }\Psi \left( \sigma \right) =\frac{1}{\Gamma \left( 1-\alpha
\right) }\int_{0}^{\sigma }\frac{\Psi ^{\prime }\left( \tau \right)}{ \left(
\sigma -\tau \right) ^{\alpha }}\,d\tau=-\Psi \left( \sigma \right)\quad(
0<\alpha <1)  \label{1.7}
\end{equation*}%
one has the solution 
\begin{equation*}
\Psi \left( \sigma \right) =E_{\alpha }\left( -\sigma ^{\alpha }\right)
\label{1.8}
\end{equation*}%
with $E_{\alpha }$ being the Mittag-Leffler function of parameter $\alpha $%
\begin{equation}
E_{\alpha }\left( z\right) =\sum_{n=0}^{\infty }\frac{z^{n}}{\Gamma \left(
\alpha n+1\right) },\quad z\in\mathbb{C}  \label{1.9}
\end{equation}%
$\left( \alpha >0\right) $. One then obtains a \textit{fractional Poisson
process\/} \cite{FCAA}, \cite{Mainardi1} with the probability of $n$ events 
\begin{equation*}
P\left( X=n\right) =\frac{\sigma ^{\alpha n}}{n!}E_{\alpha }^{(n)}\left(
-\sigma ^{\alpha }\right)  \label{1.9a}
\end{equation*}%
$E_{\alpha }^{(n)}$ denoting the $n$-th derivative of the Mittag-Leffler
function. In contrast with the Poisson case ($\alpha =1$), this process has
power law asymptotics rather than exponential, which implies that it is not
anymore Markovian. The characteristic function of this process is given by 
\begin{equation*}
C_{\alpha }\left( \lambda \right) =E_{\alpha }\left( \sigma ^{\alpha }\left(
e^{i\lambda }-1\right) \right)  \label{1.10}
\end{equation*}%
In this paper we develop an infinite-dimensional generalization of the
fractional Poisson measure and its analysis.

\section{Infinite-dimensional fractional Poisson measures}

For the Poisson measure ($\alpha =1$) an infinite-dimensional generalization
is obtained by generalizing (\ref{1.4}) to 
\begin{equation}
C\left( \varphi\right) =e^{\int \left( e^{i\varphi\left( x\right)
}-1\right)\, d\mu \left( x\right) }  \label{2.1}
\end{equation}%
for test functions $\varphi\in \mathcal{D}\left( M\right) $, $\mathcal{D}%
\left( M\right) $ being the space of $C^{\infty }$-functions of compact
support in a manifold $M$ (fixed from the very beginning), and then using
the Bochner-Minlos theorem to show that $C$ is the Fourier transform of a
measure on the distribution space $\mathcal{D}^{\prime }\left( M\right)$.
Because the Mittag-Leffler function is a \textquotedblleft
natural\textquotedblright\ generalization of the exponential function one
conjectures that an infinite-dimensional version of the fractional Poisson
measure would have a characteristic functional 
\begin{equation}
C_{\alpha }\left( \varphi\right) :=E_{\alpha }\left( \int (
e^{i\varphi\left( x\right) }-1)\,d\mu \left( x\right) \right),\quad
\varphi\in \mathcal{D}\left( M\right)  \label{2.2}
\end{equation}%
with $\mu $ a positive intensity measure fixed on the underlying manifold $M$%
. However, a priori it is not obvious that this is the Fourier transform of
a measure on $\mathcal{D}^{\prime }\left( M\right) $ nor that it corresponds
to independent processes because the Mittag-Leffler function does not
satisfy the factorization properties of the exponential.

Similarly to the Poisson case, to carry out our construction and analysis in
detail we always assume that $M$ is a geodesically complete connected
oriented (non-compact) Riemannian $C^{\infty }$-manifold, where we fix the
corresponding Borel $\sigma$-algebra $\mathcal{B}\left( M\right)$, and $\mu$
is a non-atomic Radon measure, which we assume to be non-degenerate (i.e., $%
\mu(O)>0$ for all non-empty open sets $O\subset M$). Having in mind the most
interesting applications, we also assume that $\mu(M)=\infty$.

\begin{theorem}
\label{Th3} For each $0<\alpha\leq 1$ fixed, the functional $C_{\alpha }$ in
Eq.~(\ref{2.2}) is the characteristic functional of a probability measure $%
\pi _{\mu }^{\alpha }$ on the distribution space $\mathcal{D}^{\prime
}\left( M\right) $.
\end{theorem}

\noindent \textbf{Proof.} That $C_{\alpha }$ is continuous and $C_{\alpha
}\left( 0\right) =1$ follows easily from the properties of the
Mittag-Leffler function. To check the positivity one uses the complete
monotonicity of $E_{\alpha }$, $0<\alpha <1$, which by Appendix A (Lemma \ref%
{Aug09-eq3}) implies the integral representation 
\begin{equation}
E_{\alpha }\left( -z\right) =\int_{0}^{\infty }e^{-\tau z}d\nu _{\alpha
}\left( \tau \right)  \label{2.3}
\end{equation}%
for any $z\in\mathbb{C}$ such that $\mathrm{Re}\left( z\right) \geq 0$, $%
\nu_\alpha$ being the probability measure (\ref{A2}). Hence by (\ref{2.3}) 
\begin{equation}
\sum_{a,b}C_{\alpha }\left( \varphi_{a}-\varphi_{b}\right) z_{a}^{\ast
}z_{b}=\int_{0}^{\infty }d\nu _{\alpha }\left( \tau \right)
\sum_{a,b}e^{-\tau \int_{M}d\mu \left( x\right) \left(
1-e^{i(\varphi_{a}-\varphi_{b})}\right) }z_{a}^{\ast }z_{b}  \label{2.4}
\end{equation}%
Each one of the terms in the integrand corresponds to the characteristic
function of a Poisson measure. Thus, for each $\tau$ the integrand is
positive and therefore the spectral integral (\ref{2.4}) is also positive.
From the Bochner-Minlos theorem it then follows that $C_{\alpha }$ is the
characteristic functional of a probability measure $\pi _{\mu }^{\alpha } $
on the measurable space $( \mathcal{D}^{^{\prime }}( M),\mathcal{C}_{\sigma}(%
\mathcal{D}^{^{\prime }}(M)))$, $\mathcal{C}_{\sigma}(\mathcal{D}^{^{\prime
}}(M))$ being the $\sigma$-algebra generated by the cylinder sets.

For the $\alpha=1$ case see e.g.~\cite{GV}.\hfill $\blacksquare \medskip $

Introducing the fractional Poisson measure by the above approach yields a
probability measure on $( \mathcal{D}^{^{\prime }}( M),\mathcal{C}_{\sigma}(%
\mathcal{D}^{^{\prime }}(M)))$. The next step is to find an appropriate
support for the fractional Poisson measure. Using the analyticity of the
Mittag-Leffler function one may informally rewrite (\ref{2.2}) as 
\begin{eqnarray*}
C_{\alpha }\left( \varphi\right) &=&\sum_{n=0}^{\infty }\frac{E_{\alpha
}^{\left( n\right) }\left( -\int d\mu \left( x\right) \right) }{n!}\left(
\int e^{i\varphi\left( x\right) }d\mu \left( x\right) \right) ^{n}  \notag \\
&=&\sum_{n=0}^{\infty }\frac{E_{\alpha }^{\left( n\right) }\left( -\int d\mu
\left( x\right) \right) }{n!}\int e^{i\left( \varphi\left( x_{1}\right)
+\varphi\left( x_{2}\right) +\cdots +\varphi\left( x_{n}\right) \right)
}d\mu ^{\otimes n}  \label{2.5}
\end{eqnarray*}%
For the Poisson case ($\alpha =1$) instead of $E_{\alpha }^{\left( n\right)
}\left( -\int d\mu \left( x\right) \right) $ one would have $\exp \left(
-\int d\mu \left( x\right) \right) $ for all $n$, the rest being the same.
Therefore one concludes that the main difference in the fractional case ($%
\alpha \neq 1$) is that a different weight is given to each $n$-particle
space, but that a configuration space \cite{Albeverio1}, \cite{Albeverio2}
is also the natural support of the fractional Poisson measure. The explicit
construction is made in Section 3.

Notice however that the different weights, multiplying the $n$-particle
space measures, are physically quite significant in that they have decays,
for large volumes, much smaller than the corresponding exponential factor in
the Poisson measure.

Using now the spectral representation (\ref{2.3}) of the Mittag-Leffler
function one may rewrite (\ref{2.2}) as 
\begin{equation*}
C_{\alpha }\left( \varphi \right) =\int_{0}^{\infty }\exp \left( \tau \int
(e^{i\varphi (x)}-1)\,d\mu (x)\right) \,d\nu _{\alpha }(\tau )
\end{equation*}%
with the integrand being the characteristic function of the Poisson measure $%
\pi _{\tau \mu }$, $\tau >0$. In other words, the characteristic functional (%
\ref{2.2}) coincides with the characteristic functional of the measure $%
\int_{0}^{\infty }\pi _{\tau \mu }\,d\nu _{\alpha }(\tau )$. By uniqueness,
this implies the integral decomposition 
\begin{equation*}
\pi _{\mu }^{\alpha }=\int_{0}^{\infty }\pi _{\tau \mu }\,d\nu _{\alpha
}(\tau )
\end{equation*}%
meaning that $\pi _{\mu }^{\alpha }$ is an integral (or mixture) of Poisson
measures $\pi _{\tau \mu }$, $\tau >0$.

\section{Support properties of the fractional Poisson measure}

\subsection{Configuration spaces}

\bigskip The \textit{configuration space} $\Gamma:=\Gamma _{M}$ over the
manifold $M$ is defined as the set of all locally finite subsets of $M$
(simple configurations), 
\begin{equation}
\Gamma:=\{\gamma \subset M\,:\,|\gamma \cap K|<\infty \;\mathrm{%
for\;any\;compact\;}K\subset M\}  \label{2eq70}
\end{equation}%
Here (and below) $|A|$ denotes the cardinality of a set $A$.

As usual we identify each $\gamma \in \Gamma $ with a non-negative
integer-valued Radon measure, 
\begin{equation*}
\Gamma \ni \gamma \mapsto \sum_{x\in \gamma }\delta _{x}\in \mathcal{M}(M)
\end{equation*}%
where $\delta _{x}$ is the Dirac measure with unit mass at $x$, $\sum_{x\in
\emptyset }\delta _{x}:=$ zero measure, and $\mathcal{M}(M)$ denotes the set
of all non-negative Radon measures on $\mathcal{B}\left( M\right) $. In this
way the space $\Gamma $ can be endowed with the relative topology as a
subset of the space $\mathcal{M}(M)$ with the vague topology, i.e., the
weakest topology on $\Gamma $ for which the mappings 
\begin{equation*}
\Gamma \ni \gamma \mapsto \left\langle \gamma ,f\right\rangle
:=\int_{M}f(x)d\gamma (x)=\sum_{x\in \gamma }f(x)
\end{equation*}%
are continuous for all real-valued continuous functions $f$ on $M$ with
compact support. We denote the corresponding Borel $\sigma $-algebra on $%
\Gamma $ by $\mathcal{B}(\Gamma )$ .

For each $Y\in \mathcal{B}(M)$ let us consider the space $\Gamma _{Y}$ of
all configurations contained in $Y$, 
\begin{equation*}
\Gamma _{Y}:=\left\{ \gamma \in \Gamma :\left\vert \gamma \cap (X\backslash
Y)\right\vert =0\right\}
\end{equation*}%
and the space $\Gamma _{Y}^{(n)}$ of $n$-point configurations, 
\begin{equation*}
\Gamma _{Y}^{(n)}:=\left\{ \gamma \in \Gamma _{Y}:\left\vert \gamma
\right\vert =n\right\} ,n\in \mathbb{N},\quad \Gamma _{Y}^{(0)}:=\left\{
\emptyset \right\}
\end{equation*}%
A topological structure may be introduced on $\Gamma _{Y}^{(n)}$ through the
natural surjective mapping of 
\begin{equation*}
\widetilde{Y^{n}}:=\left\{ (x_{1},...,x_{n}):x_{i}\in Y,x_{i}\neq x_{j}%
\hbox{ if }i\neq j\right\}
\end{equation*}%
onto $\Gamma _{Y}^{(n)}$, 
\begin{equation*}
\begin{array}{ll}
\mathrm{sym}_{Y}^{n}:\widetilde{Y^{n}} & \longrightarrow \Gamma _{Y}^{(n)}
\\ 
(x_{1},...,x_{n}) & \longmapsto \left\{ x_{1},...,x_{n}\right\}%
\end{array}%
\end{equation*}%
which is at the origin of a bijection between $\Gamma _{Y}^{(n)}$ and the
symmetrization $\widetilde{Y^{n}}\diagup S_{n}$ of $\widetilde{Y^{n}}$, $%
S_{n}$ being the permutation group over $\left\{ 1,...,n\right\} $. Thus, $%
\mathrm{sym}_{Y}^{n}$ induces a metric on $\Gamma _{Y}^{(n)}$ and the
corresponding Borel $\sigma $-algebra $\mathcal{B}(\Gamma _{Y}^{(n)})$ on $%
\Gamma _{Y}^{(n)}$.

For $\Lambda \in \mathcal{B}(M)$ with compact closure ($\Lambda \in \mathcal{%
B}_c(M)$ for short), it clearly follows from (\ref{2eq70}) that 
\begin{equation*}
\Gamma _\Lambda =\bigsqcup_{n=0}^\infty \Gamma _\Lambda ^{(n)}
\end{equation*}
One defines the $\sigma $-algebra $\mathcal{B}(\Gamma _\Lambda )$ by the
disjoint union of the $\sigma $-algebras $\mathcal{B}(\Gamma _\Lambda
^{(n)}) $, $n\in \mathbb{N}_0 $.

For each $\Lambda \in \mathcal{B}_{c}(M)$ there is a natural measurable
mapping $p_{\Lambda }:\Gamma \rightarrow \Gamma _{\Lambda }$. Similarly,
given any pair $\Lambda _{1},\Lambda _{2}\in \mathcal{B}_{c}(M)$ with $%
\Lambda _{1}\subset \Lambda _{2}$ there is a natural mapping $p_{\Lambda
_{2},\Lambda _{1}}:\Gamma _{\Lambda _{2}}\rightarrow \Gamma _{\Lambda _{1}}$%
. They are defined, respectively, by 
\begin{equation*}
\begin{array}{ll}
p_{\Lambda }: & \Gamma \longrightarrow \Gamma _{\Lambda } \\ 
& \gamma \longmapsto \gamma _{\Lambda }:=\gamma \cap \Lambda%
\end{array}%
\quad 
\begin{array}{ll}
p_{\Lambda _{2},\Lambda _{1}}: & \Gamma _{\Lambda _{2}}\longrightarrow
\Gamma _{\Lambda _{1}} \\ 
& \gamma \longmapsto \gamma _{\Lambda _{1}}%
\end{array}%
\end{equation*}%
It can be shown (cf.~\cite{S94}) that $(\Gamma ,\mathcal{B}(\Gamma ))$
coincides (up to an isomorphism) with the projective limit of the measurable
spaces $(\Gamma _{\Lambda },\mathcal{B}(\Gamma _{\Lambda }))$, $\Lambda \in 
\mathcal{B}_{c}(M)$, with respect to the projection $p_{\Lambda }$, i.e., $%
\mathcal{B}(\Gamma )$ is the smallest $\sigma $-algebra on $\Gamma $ with
respect to which all projections $p_{\Lambda }$, $\Lambda \in \mathcal{B}%
_{c}(M)$, are measurable.

\subsection{Fractional Poisson measure on $\Gamma$}

Given a measure $\mu $ on the underlying measurable space $(M,\mathcal{B}%
(M)) $ described before, consider for each $n\in \mathbb{N}$ the product
measure $\mu ^{\otimes n}$ on $(M^{n},\mathcal{B}(M^{n}))$. Since $\mu
^{\otimes n}(M^{n}\backslash \widetilde{M^{n}})=0$, one may consider for
each $\Lambda \in \mathcal{B}_{c}(M)$ the restriction of $\mu ^{\otimes }$
to $(\widetilde{\Lambda ^{n}},\mathcal{B}(\widetilde{\Lambda ^{n}}))$, which
is a finite measure, and then the image measure $\mu _{\Lambda }^{(n)}$ on $%
(\Gamma _{\Lambda }^{(n)},\mathcal{B}(\Gamma _{\Lambda }^{(n)}))$ under the
mapping $\mathrm{sym}_{\Lambda }^{n}$, 
\begin{equation*}
\mu _{\Lambda }^{(n)}:=\mu ^{\otimes n}\circ (\mathrm{sym}_{\Lambda
}^{n})^{-1}
\end{equation*}%
For $n=0$ we set $\mu _{\Lambda }^{(0)}:=1$\footnote{%
Of course this construction holds for any Borel set $Y\in \mathcal{B}(M)$.
In this case, $\mu _{Y}^{(n)}(\Gamma _{Y}^{(n)})<\infty $ provided $\mu
(Y)<\infty $. For more details and proofs see e.g.~\cite{YuKu}, \cite%
{KoKu100}.}. Now, for each $0<\alpha <1$ one may define a probability
measure $\pi _{\mu ,\Lambda }^{\alpha }$ on $(\Gamma _{\Lambda },\mathcal{B}%
(\Gamma _{\Lambda }))$ by 
\begin{equation}
\pi _{\mu ,\Lambda }^{\alpha }:=\sum_{n=0}^{\infty }\frac{E_{\alpha
}^{(n)}(-\mu (\Lambda ))}{n!}\mu _{\Lambda }^{(n)}  \label{t94}
\end{equation}

The family $\{\pi _{\mu ,\Lambda }^{\alpha }:\Lambda \in \mathcal{B}%
_{c}(M)\} $ of probability measures yields a probability measure on $(\Gamma
,\mathcal{B}(\Gamma ))$. In fact, this family is consistent, that is, 
\begin{equation*}
\pi _{\mu ,\Lambda _{1}}^{\alpha }=\pi _{\mu ,\Lambda _{2}}^{\alpha }\circ
p_{\Lambda _{2},\Lambda _{1}}^{-1},\quad \forall \,\Lambda _{1},\Lambda
_{2}\in \mathcal{B}_{c}(M),\Lambda _{1}\subset \Lambda _{2}
\end{equation*}%
and thus, by the version of Kolmogorov's theorem for the projective limit
space $(\Gamma ,\mathcal{B}(\Gamma ))$ \cite[Chap.~V Theorem~5.1]{P67}, the
family $\{\pi _{\mu ,\Lambda }^{\alpha }:\Lambda \in \mathcal{B}_{c}(M)\}$
determines uniquely a measure $\pi _{\mu }^{\alpha }$ on $(\Gamma ,\mathcal{B%
}(\Gamma ))$ such that 
\begin{equation*}
\pi _{\mu ,\Lambda }^{\alpha }=\pi _{\mu }^{\alpha }\circ p_{\Lambda
}^{-1},\quad \forall \,\Lambda \in \mathcal{B}_{c}(M)
\end{equation*}

Let us now compute the characteristic functional of the measure $\pi _{\mu
}^{\alpha}$. Given a $\varphi\in \mathcal{D}(M)$ we have supp$%
\,\varphi\subset \Lambda $ for some $\Lambda \in \mathcal{B}_{c}(M)$,
meaning that 
\begin{equation*}
\langle \gamma ,\varphi\rangle =\langle p_\Lambda(\gamma),\varphi\rangle ,
\quad\forall\,\gamma\in \Gamma
\end{equation*}%
Thus 
\begin{equation*}
\int_{\Gamma }e^{i\langle \gamma ,\varphi\rangle }d\pi _{\mu }^{\alpha
}(\gamma )=\int_{\Gamma _{\Lambda }}e^{i\langle \gamma ,\varphi\rangle }d\pi
_{\mu,\Lambda}^{\alpha }(\gamma )
\end{equation*}%
and the infinite divisibility (\ref{t94}) of the measure $\pi _{\mu,\Lambda
}^{\alpha}$ yields for the right-hand side of the equality 
\begin{equation*}
\sum_{n=0}^{\infty }\frac{E_{\alpha }^{(n)}(-\mu (\Lambda ))}{n!}%
\int_{\Lambda ^{n}}e^{i(\varphi(x_{1})+\ldots +\varphi(x_{n}))}d\mu
^{\otimes n}(x)=\sum_{n=0}^{\infty }\frac{E_{\alpha }^{(n)}(-\mu (\Lambda ))%
}{n!}\left( \int_{\Lambda }e^{i\varphi(x)}d\mu (x)\right) ^{n}
\end{equation*}%
which corresponds to the Taylor expansion of the function 
\begin{equation*}
E_{\alpha }\left( \int_{\Lambda }(e^{i\varphi(x)}-1)\,d\mu (x)\right)=
E_{\alpha }\left( \int_{M}(e^{i\varphi(x)}-1)\,d\mu (x)\right)
\end{equation*}%
In other words, the characteristic functional of the measure $\pi _{\mu
}^{\alpha }$ coincides with the characteristic functional of the probability
measure given by Theorem \ref{Th3} through the Bochner-Minlos theorem.

Similarly to the $\alpha=1$ case, this shows that the probability measure on 
$( \mathcal{D}^{^{\prime }}( M),\mathcal{C}_{\sigma}( \mathcal{D}^{^{\prime
}}( M))) $ given by Theorem \ref{Th3} is actually supported on generalized
functions of the form $\sum_{x\in\gamma}\delta_x$, $\gamma\in\Gamma$. Thus,
each fractional Poisson measure $\pi _{\mu }^{\alpha }$ can either be
consider on $(\Gamma ,\mathcal{B}(\Gamma ))$ or on $(\mathcal{D}^{\prime },%
\mathcal{C}_{\sigma}(\mathcal{D}^{\prime }(M)))$ where, in contrast to $%
\Gamma $, $\mathcal{D}^{\prime }(M)\supset \Gamma$ is a linear space. Since $%
\pi _{\mu }^{\alpha }(\Gamma )=1$, the measure space $(\mathcal{D}^{\prime
}(M),\mathcal{C}_{\sigma}(\mathcal{D}^{\prime }(M)),\pi _{\mu }^{\alpha })$
can, in this way, be regarded as a linear extension of the fractional
Poisson space $(\Gamma ,\mathcal{B}(\Gamma ),\pi _{\mu }^{\alpha })$.

\section{Fractional Poisson analysis\label{Section4}}

\subsection{Fractional Lebesgue-Poisson measure and unitary isomorphisms}

Let us now consider the space of finite configurations 
\begin{equation*}
\Gamma _{0}:=\bigsqcup_{n=0}^{\infty }\Gamma _{M}^{(n)}
\end{equation*}%
endowed with the topology of disjoint union of topological spaces, with the
corresponding Borel $\sigma $-algebra $\mathcal{B}(\Gamma _{0})$ and the
so-called $K$-transform \cite{YuKu}, \cite{KKO04}, \cite{Le72}, \cite{Le75a}%
, \cite{Le75b}, a mapping which maps functions defined on $\Gamma _{0}$ into
functions defined on $\Gamma $. By definition, given a $\mathcal{B}(\Gamma
_{0})$-measurable function $G$ with local support, that is, $%
G\!\!\upharpoonright _{\Gamma _{0}\backslash \Gamma _{\Lambda }}\equiv 0$
for some $\Lambda \in \mathcal{B}_{c}(M)$, the $K$-transform of $G$ is a
mapping $KG:\Gamma \rightarrow \mathbb{R}$ defined at each $\gamma \in
\Gamma $ by 
\begin{equation}
(KG)(\gamma ):=\sum_{{\eta \subset \gamma } \atop {|\eta |<\infty }}%
G(\eta )  \label{Eq2.9}
\end{equation}%
Note that for every such function $G$ the sum in (\ref{Eq2.9}) has only a
finite number of summands different from zero, and thus $KG$ is a
well-defined function on $\Gamma $. Moreover, if $G$ has support described
as before, then the restriction $(KG)\!\!\upharpoonright _{\Gamma _{\Lambda
}}$ is a $\mathcal{B}(\Gamma _{\Lambda })$-measurable function and $%
(KG)(\gamma )=(KG)\!\!\upharpoonright _{\Gamma _{\Lambda }}\!\!(\gamma
_{\Lambda })$ for all $\gamma \in \Gamma $.

In terms of the dual operator $K^{\ast }$ of the $K$-transform, this means
that the image of a probability measure on $\Gamma $ under $K^{\ast }$
yields a measure on $\Gamma _{0}$. More precisely, given a probability
measure $\nu $ on $(\Gamma ,\mathcal{B}(\Gamma ))$ with finite local moments
of all orders, that is, 
\begin{equation*}
\int_{\Gamma }|\gamma _{\Lambda }|^{n}\,d\nu (\gamma )<\infty \quad \mathrm{%
for\,\,all}\,\,n\in \mathbb{N}\mathrm{\,\,and\,\,all\,\,}\Lambda \in 
\mathcal{B}_{c}(M)
\end{equation*}%
then $K^{\ast }\nu $ is a measure on $(\Gamma _{0},\mathcal{B}(\Gamma _{0}))$
defined on each bounded $\mathcal{B}(\Gamma _{0})$-measurable set $A$ by 
\begin{equation*}
(K^{\ast }\nu )(A)=\int_{\Gamma }(K\boldsymbol{1}_{A})(\gamma )\,d\nu
(\gamma )
\end{equation*}%
The measure $K^{\ast }\nu $ is called the correlation measure corresponding
to $\nu $. In particular, for the Poisson measure $\pi _{\mu }$, the
correlation measure corresponding to $\pi _{\mu }$ is called the
Lebesgue-Poisson measure 
\begin{equation*}
\lambda _{\mu }:=\sum_{n=0}^{\infty }\frac{1}{n!}\mu ^{(n)},\quad \mu
^{(n)}:=\mu ^{\otimes n}\circ (\mathrm{sym}_{M}^{n})^{-1}
\end{equation*}%
For more details and proofs see e.g.~\cite{YuKu}.

\begin{theorem}
\label{Th2}For each $0<\alpha <1$, the correlation measure corresponding to
the fractional Poisson measure $\pi _{\mu }^{\alpha }$ is the measure on $%
(\Gamma _{0},\mathcal{B}(\Gamma _{0}))$ given by 
\begin{equation}
\lambda _{\mu }^{\alpha }:=\sum_{n=0}^{\infty }\frac{1}{\Gamma (\alpha n+1)}%
\mu ^{(n)}  \label{forma}
\end{equation}%
In other words, $d\lambda _{\mu }^{\alpha }=E_{\alpha }^{(|\cdot
|)}(0)\,d\lambda _{\mu }$.
\end{theorem}

In the sequel we call the measure $\lambda_\mu^\alpha$ the \textit{%
fractional Lebesgue-Poisson measure\/}.

\medskip

\noindent \textbf{Proof.} Let $A$ be a bounded $\mathcal{B}(\Gamma _{0})$%
-measurable set, that is, 
\begin{equation*}
A\subset \bigsqcup_{n=0}^{N}\Gamma _{\Lambda }^{(n)}
\end{equation*}%
for some $N\in \mathbb{N}_{0}$ and some $\Lambda \in \mathcal{B}_{c}(M)$. By
the previous considerations, this means that for all $\gamma \in \Gamma $
one has $(K\boldsymbol{1}_{A})(\gamma )=(K\boldsymbol{1}_{A})(\gamma
_{\Lambda })$, and thus 
\begin{eqnarray*}
\int_{\Gamma }(K\boldsymbol{1}_{A})(\gamma )\,d\pi _{\mu }^{\alpha }(\gamma
) &=&\int_{\Gamma _{\Lambda }}(K\boldsymbol{1}_{A})(\gamma )\,d\pi _{\mu
,\Lambda }^{\alpha }(\gamma ) \\
&=&\sum_{n=0}^{\infty }\frac{E_{\alpha }^{(n)}(-\mu (\Lambda ))}{n!}%
\int_{\Gamma _{\Lambda }^{(n)}}(K\boldsymbol{1}_{A})(\eta )\,d\mu _{\Lambda
}^{(n)}\left( \eta \right) \\
&=&\int_{\Gamma _{\Lambda }}E_{\alpha }^{(|\eta |)}(-\mu (\Lambda ))(K%
\boldsymbol{1}_{A})(\eta )\,d\lambda _{\mu }(\eta )
\end{eqnarray*}%
Observe that the latter integral is with respect to the Lebesgue-Poisson
measure $\lambda _{\mu }$, which properties are well-known (see e.g.~\cite%
{YuKu}). In particular, those yield 
\begin{eqnarray*}
&&\int_{\Gamma _{\Lambda }}E_{\alpha }^{(|\eta |)}(-\mu (\Lambda ))(K%
\boldsymbol{1}_{A})(\eta )\,d\lambda _{\mu }(\eta ) \\
&=&\int_{\Gamma _{0}}E_{\alpha }^{(|\eta |)}(-\mu (\Lambda ))\sum_{\xi
\subset \eta }\boldsymbol{1}_{A}(\xi )\boldsymbol{1}_{\Gamma _{\Lambda
}}(\eta \setminus \xi )\,d\lambda _{\mu }(\eta ) \\
&=&\int_{\Gamma _{0}}\boldsymbol{1}_{A}(\eta )\left( \int_{\Gamma
_{0}}E_{\alpha }^{(|\eta \cup \xi |)}(-\mu (\Lambda ))\boldsymbol{1}_{\Gamma
_{\Lambda }}(\xi )\,d\lambda _{\mu }(\xi )\right) d\lambda _{\mu }(\eta )
\end{eqnarray*}%
where for each $\eta \in \Gamma _{0}$ fixed, i.e., $\eta \in \Gamma
_{M}^{(m)}$ for some $m\in \mathbb{N}_{0}$, the integral between brackets is
given by 
\begin{eqnarray*}
&&\sum_{n=0}^{\infty }\frac{1}{n!}\int_{\Gamma _{M}^{(n)}}E_{\alpha
}^{(|\eta \cup \xi |)}(-\mu (\Lambda ))\boldsymbol{1}_{\Gamma _{\Lambda
}}(\xi )\,d\mu ^{(n)}(\xi ) \\
&=&\sum_{n=0}^{\infty }\frac{E_{\alpha }^{(m+n)}(-\mu (\Lambda ))}{n!}(\mu
(\Lambda ))^{n} \\
&=&E_{\alpha }^{(m)}(-\mu (\Lambda )+\mu (\Lambda ))
\end{eqnarray*}%
As a result, 
\begin{equation*}
\int_{\Gamma }(K\boldsymbol{1}_{A})(\gamma )\,d\pi _{\mu }^{\alpha }(\gamma
)=\int_{\Gamma _{0}}\boldsymbol{1}_{A}(\eta )E_{\alpha }^{(|\eta
|)}(0)\,d\lambda _{\mu }(\eta )
\end{equation*}%
showing that the correlation measure corresponding to $\pi _{\mu }^{\alpha }$
is absolutely continuous with respect to the Lebesgue-Poisson measure $%
\lambda _{\mu }$. Moreover, denoting such a correlation measure by $\lambda
_{\mu }^{\alpha }$, the density is given by $\frac{d\lambda _{\mu }^{\alpha }%
}{d\lambda _{\mu }}=E_{\alpha }^{(|\cdot |)}(0)$.

To conclude, notice that for each $n\in \mathbb{N}_{0}$ one has 
\begin{equation*}
E_{\alpha }^{(n)}(0)=\frac{n!}{\Gamma (\alpha n+1)}
\end{equation*}%
which combined with the definition of the measure $\lambda _{\mu }$ leads to
(\ref{forma}).\hfill $\blacksquare \medskip $

Throughout this work all $L^{p}$-spaces consist of complex-valued functions.
For simplicity, the $L^{p}$-spaces with respect to a measure $\nu $ will be
shortly denoted by $L^{p}(\nu )$ if the underlying measurable space is clear
from the context.

\begin{corollary}
We have $G\in L^{2}(\lambda _{\mu }^{\alpha })$ if and only if $G\sqrt{%
E_{\alpha }^{(|\cdot |)}(0)}\in L^{2}(\lambda _{\mu })$. Then, 
\begin{equation*}
\Vert G\Vert _{L^{2}(\lambda _{\mu }^{\alpha })}=\left\Vert G\sqrt{E_{\alpha
}^{(|\cdot |)}(0)}\right\Vert _{L^{2}(\lambda _{\mu })}
\end{equation*}
\end{corollary}

This result states that there is a unitary isomorphism between the spaces $%
L^{2}(\lambda _{\mu }^{\alpha })$ and $L^{2}(\lambda _{\mu })$: 
\begin{eqnarray*}
I_{\alpha }:L^{2}(\lambda _{\mu }^{\alpha }) &\rightarrow &L^{2}(\lambda
_{\mu }) \\
I_{\alpha }(G):= &&G\sqrt{E_{\alpha }^{(|\cdot |)}(0)}
\end{eqnarray*}%
Hence, through $I_{\alpha }$ one may extend the unitary isomorphisms defined
between the space $L^{2}(\lambda _{\mu })$ and the (Bose or symmetric) Fock
space $\mathrm{Exp}L^{2}(\mu )$ and between the space $L^{2}(\lambda _{\mu
}) $ and $L^{2}(\pi _{\mu })$ \cite{KoKu100}, \cite{Oliveira1} to $%
L^{2}(\lambda _{\mu }^{\alpha })$, $0<\alpha \leq 1$: 
\begin{equation*}
\begin{matrix}
L^{2}(\lambda _{\mu }^{\alpha }) & \overset{I_{\alpha }}{\mapsto } & 
L^{2}(\lambda _{\mu }) & \overset{I_{\lambda \pi }}{\mapsto } & L^{2}(\pi
_{\mu }) & \overset{I_{\pi }}{\mapsto } & \mathrm{Exp}L^{2}(\mu ) \\ 
G & \mapsto & G\sqrt{E_{\alpha }^{(|\cdot |)}(0)} & \mapsto & 
\sum_{n=0}^{\infty }\langle C_{n}^{\mu },g^{(n)}\rangle & \mapsto & \left(
g^{(n)}\right) _{n=0}^{\infty }%
\end{matrix}%
\end{equation*}%
for 
\begin{equation*}
g^{(n)}(x_{1},\ldots ,x_{n}):=\frac{\sqrt{E_{\alpha }^{(n)}(0)}}{n!}%
G(\{x_{1},\ldots ,x_{n}\}),\quad g^{(0)}:=E_{\alpha }(0)G(\emptyset
)=G(\emptyset )
\end{equation*}%
and $C_{n}^{\mu }$ a Charlier kernel.

In particular, the image of a Fock coherent state $e(f):=(\frac{f^{\otimes n}%
}{n!})_{n=0}^\infty$, $f\in L^2(\mu)$, under $(I_\pi\circ
I_{\lambda\pi})^{-1}$ is the (Lebesgue-Poisson) coherent state $%
e_\lambda(f):\Gamma _0\rightarrow \mathbb{C}$ defined for any $\mathcal{B}%
(M) $-measurable function $f:M\to\mathbb{C}$ by 
\begin{equation*}
e_\lambda(f,\eta ):=\prod_{x\in \eta }f\left( x\right) ,\ \eta \in \Gamma
_0\setminus\{\emptyset\},\quad e_\lambda(f,\emptyset):=1
\end{equation*}
This definition implies that $e_\lambda(f)\in L^p(\lambda_\mu)$ whenever $%
f\in L^p(\mu)$ for some $p\geq 1$. Moreover, $\|e_\lambda(f)\|^p_{L^p(%
\lambda_\mu)}=\exp\left(\|f\|^p_{L^p(\mu)}\right)$. For $\alpha\not=1$, the
following result holds.

\begin{proposition}
Let $0<\alpha<1$ and $p\geq 1$ be given. For all $f\in L^p(\mu)$ we have $%
e_\lambda(f)\in L^p(\lambda_\mu^\alpha)$ and 
\begin{equation*}
\|e_\lambda(f)\|^p_{L^p(\lambda_\mu^\alpha)}=E_\alpha\left(\|f\|^p_{L^p(%
\mu)}\right)
\end{equation*}
\end{proposition}

\noindent \textbf{Proof.} By Theorem \ref{Th2}, given a $f\in L^{p}(\mu )$
for some $p\geq 1$, 
\begin{eqnarray*}
\Vert e_{\lambda }(f)\Vert _{L^{p}(\lambda _{\mu }^{\alpha })}^{p}
&=&\int_{\Gamma _{0}}|e_{\lambda }(f,\eta )|^{p}E_{\alpha }^{(|\eta
|)}(0)\,d\lambda _{\mu }(\eta ) \\
&=&\sum_{n=0}^{\infty }\frac{1}{\Gamma (\alpha n+1)}\left(
\int_{M}|f(x)|^{p}\,d\mu (x)\right) ^{n}
\end{eqnarray*}%
which by the Taylor expansion (\ref{1.9}) is equal to $E_{\alpha }\left(
\int_{M}|f(x)|^{p}\,d\mu (x)\right).\hfill\blacksquare \medskip $

According to the latter considerations, the realization of a coherent state $%
e(f)$, $f\in L^{2}(\mu )$, in a $L^{2}(\lambda _{\mu }^{\alpha })$ space is $%
\lambda_{\mu }$-a.e.~given by 
\begin{equation}
I_{\alpha }^{-1}e_{\lambda }(f)=\frac{e_{\lambda }(f)}{\sqrt{E_{\alpha
}^{(|\cdot |)}(0)}}  \label{qqcoisa}
\end{equation}%
In addition, given a dense subspace $L\subseteq L^2(\mu)$, the set $%
\{I_{\alpha }^{-1}e_{\lambda }(f): f\in L\}$ is total in $L^{2}(\lambda
_{\mu }^{\alpha })$. As in the Lebesgue-Poisson case, we define the
fractional (Lebesgue-Poisson) coherent state $e_{\alpha }(f):\Gamma
_{0}\rightarrow \mathbb{C}$ corresponding to a $\mathcal{B}(M)$-measurable
function $f$ by 
\begin{equation*}
e_{\alpha }(f,\eta ):=\frac{e_{\lambda }(f,\eta)}{\sqrt{E_{\alpha }^{(|\eta
|)}(0)}},\quad \forall \,\eta \in \Gamma _{0}
\end{equation*}

\subsection{Annihilation and creation operators}

The unitary isomorphism between the Fock space and $L^{2}(\lambda _{\mu })$
provides natural operators on the space $L^{2}(\lambda _{\mu })$ by carrying
over the standard Fock space operators. In particular, the annihilation and
the creation operators, for which the images in $L^{2}(\lambda _{\mu })$ are
well-known, see e.g.~\cite{FicWin93}, 
\begin{equation*}
\left( a_{\lambda }^{-}(\varphi )G\right) (\eta ):=\int_{M}G(\eta \cup
\{x\})\varphi (x)\,d\mu (x),\quad \eta \in \Gamma _{0}
\end{equation*}%
and 
\begin{equation*}
\left( a_{\lambda }^{+}(\varphi )G\right) (\eta ):=\sum_{x\in \eta }G(\eta
\backslash \{x\})\varphi (x),\quad \lambda _{\mu }-a.a.\,\eta \in \Gamma _{0}
\end{equation*}%
Here $\varphi \in \mathcal{D}(M)$ and $G$ is a complex-valued bounded $%
\mathcal{B}(\Gamma _{0})$-measurable function with bounded support, i.e., $%
G\!\!\upharpoonright _{\Gamma _{0}\backslash \left(
\bigsqcup_{n=0}^{N}\Gamma _{\Lambda }^{(n)}\right) }\equiv 0$ for some $%
\Lambda \in \mathcal{B}_{c}(M)$ and some $N\in \mathbb{N}_{0}$. In the
sequel we denote the space of such functions $G$ by $B_{bs}(\Gamma _{0})$.

For more details and proofs see e.g.~\cite{KoKu100}, \cite{Oliveira1} and
the references therein.

Through the unitary isomorphism $I_{\alpha }^{-1}$, $0<\alpha <1$, the same
Fock space operators can naturally be carried over to the space $%
L^{2}(\lambda _{\mu }^{\alpha })$.

\begin{proposition}
\label{ancr} For each $\varphi \in \mathcal{D}(M)$, the following relations
hold on $B_{bs}(\Gamma _{0})$: 
\begin{equation*}
a^-_\alpha(\varphi):=I_\alpha^{-1}a^-_\lambda(\varphi)I_\alpha =\sqrt{\frac{%
E_\alpha^{(|\cdot|+1)}(0)}{E_\alpha^{(|\cdot|)}(0)}}a^-_\lambda(\varphi)
\end{equation*}
and 
\begin{equation*}
a^+_\alpha(\varphi):=I_\alpha^{-1}a^+_\lambda(\varphi)I_\alpha =\sqrt{\frac{%
E_\alpha^{(|\cdot|-1)}(0)}{E_\alpha^{(|\cdot|)}(0)}}a^+_\lambda(\varphi)
\end{equation*}
\end{proposition}

\noindent \textbf{Proof.} One first observes that $I_{\alpha }$ maps the
space $B_{bs}(\Gamma _{0})$ into itself. In fact, given a $G\in
B_{bs}(\Gamma _{0})$, i.e., $G\!\!\upharpoonright _{\Gamma _{0}\backslash
\left( \bigsqcup_{n=0}^{N}\Gamma _{\Lambda }^{(n)}\right) }\equiv 0$ for
some $\Lambda \in \mathcal{B}_{c}(M)$ and some $N\in \mathbb{N}_{0}$, one
has 
\begin{equation*}
|(I_{\alpha }G)(\eta )|=\sqrt{E_{\alpha }^{(|\eta |)}(0)}|G(\eta )|\leq
\max_{0\leq n\leq N}\sqrt{\frac{n!}{\Gamma (\alpha n+1)}}\sup_{\eta \in
\Gamma _{0}}(|G(\eta )|),\quad \forall \,\eta \in \Gamma _{0}
\end{equation*}%
showing that $I_{\alpha }G$ is bounded. Since the support of $I_{\alpha }G$
clearly coincides with the support of $G$, this means that $I_{\alpha }G\in
B_{bs}(\Gamma _{0})$.

Hence, given a $G\in B_{bs}(\Gamma_0)$, for all $\eta\in\Gamma_0$ one has 
\begin{eqnarray*}
(a^-_\lambda(\varphi)(I_\alpha G))(\eta)&=&\int_M (I_\alpha G)(\eta \cup
\{x\})\varphi (x)\,d\mu(x) \\
&=&\sqrt{E_\alpha^{(|\eta|+1)}(0)}\,(a^-_\lambda(\varphi)G)(\eta)
\end{eqnarray*}
which proves the first equality by calculating the image of both sides under 
$I_\alpha^{-1}$. A similar procedure applied to $a^+_\alpha(\varphi)$
completes the proof.\hfill$\blacksquare \medskip$

\subsection{Second quantization operators}

Given a contraction operator $B$ on $L^{2}(\mu )$ one may define a
contraction operator $\mathrm{Exp}B$ on the Fock space $\mathrm{Exp}%
L^{2}(\mu )$ acting on coherent states $e(f)$, $f\in L^{2}(\mu )$, by $%
\mathrm{Exp}B\left( e(f)\right) =e(Bf)$. In particular, given a positive
self-adjoint operator $A$ on $L^{2}(\mu )$ and the contraction semigroup $%
e^{-tA}$, $t\geq 0$, one can define a contraction semigroup $\mathrm{Exp}%
\left( e^{-tA}\right) $ on $\mathrm{Exp}L^{2}(\mu )$ in this way. The
generator is the well-known second quantization operator corresponding to $A$%
. We denote it by $d\mathrm{Exp}A$. Through the unitary isomorphism between
the Fock space and the space $L^{2}(\lambda _{\mu })$ one may then define
the corresponding operator in $L^{2}(\lambda _{\mu })$. We denote the
(Lebesgue-Poisson) second quantization operator corresponding to $A$ by $%
H_{A}^{LP}$. The action of $H_{A}^{LP}$ on coherent states is given by 
\begin{equation*}
\left( H_{A}^{LP}e_{\lambda }(f)\right) (\eta )=\sum_{x\in \eta }\left(
Af\right) (x)e_{\lambda }(f,\eta \backslash \{x\}),\quad f\in D(A)
\end{equation*}

Through the unitary isomorphism $I_{\alpha }^{-1}$, $0<\alpha <1$, the
second quantization operator can also be carried over to the space $%
L^{2}(\lambda _{\mu }^{\alpha })$: 
\begin{equation*}
H_{A}^{\alpha }:=I_{\alpha }^{-1}H_{A}^{LP}I_{\alpha }
\end{equation*}

\begin{proposition}
For any $f\in D(A)$ we have 
\begin{equation*}
\left( H_{A}^{\alpha }e_{\alpha }(f)\right) (\eta )=\sqrt{\frac{E_{\alpha
}^{(|\eta |-1)}(0)}{E_{\alpha }^{(|\eta |)}(0)}}\sum_{x\in \eta }\left(
Af\right) (x)e_{\alpha }(f,\eta \backslash \{x\})
\end{equation*}
\end{proposition}

\noindent \textbf{Proof.} According to (\ref{qqcoisa}), $I_\alpha
e_\alpha(f)=e_\lambda(f)$, and thus for $\lambda_\mu$-a.a.~$\eta\in\Gamma_0$%
, 
\begin{equation*}
\left(H_A^{LP}(I_\alpha e_\alpha(f))\right)(\eta) =\left(
H_{A}^{LP}e_{\lambda }(f)\right) (\eta ) =\sqrt{E_{\alpha}^{(|\eta |-1)}(0)}%
\sum_{x\in \eta }\left( Af\right) (x)e_{\alpha}(f,\eta \backslash \{x\})
\end{equation*}
leading to the required result by calculating the image of both sides under $%
I_\alpha^{-1}$.\hfill$\blacksquare \medskip$

\section{Conclusions}

Replacing the exponential, in the characteristic functional (\ref{2.1}) of
the infinite-dimensional Poisson measure, by a Mittag-Leffler function one
obtains the characteristic functional of a consistent measure in the
distribution space $\mathcal{D}^{\prime }\left( M\right) $. As for the
infinite-dimensional Poisson measure the support of this new measure is
spanned by distributions of the form $\sum \delta _{x}$, implying that it
may also be interpreted as a measure in configuration spaces.

The identity of the supports allows for the development of a fractional
infinite-dimensional analysis modeled on Poisson analysis. Although the
support is the same, the new measure displays some noticeable differences in
relation to the Poisson measure, namely, the much slower rate of decay of
the weights for the $n$-particle space measures. This might have physical
consequences, for example when such measures are used to describe
interacting particle systems.

The different weight $E_{\alpha }^{(n)}(-\mu (\Lambda ))$, given to each $n-$%
particle space, as opposed to the uniform $\exp (-\mu (\Lambda ))$ of the
Poisson case, also implies that through the isomorphism of Section 4 one
obtains an interacting Fock space.

\section*{Appendix A. Complete monotonicity of the Mittag-Leffler function
for complex arguments}

A positive $C^{\infty }$-function $f$ is said to be completely monotone if
for each $k\in \mathbb{N}_{0}$ 
\begin{equation*}
(-1)^{k}f^{(k)}(t)\geq 0,\quad \forall t>0
\end{equation*}%
According to Bernstein's theorem (see e.g.~\cite[Chapter XIII.4 Theorem~1]%
{Fel71}), for functions $f$ such that $f(0^{+})=1$ the complete monotonicity
property is equivalent to the existence of a probability measure $\nu $ on $%
\mathbb{R}_{0}^{+}$ such that 
\begin{equation*}
f(t)=\int_{0}^{\infty }e^{-t\tau }\,d\nu (\tau )<\infty ,\quad \forall \,t>0
\end{equation*}%
H.~Pollard in \cite{Pollard} proved the complete monotonicity of $E_{\alpha
} $, $0<\alpha <1$, for non-positive real arguments showing that 
\begin{equation}
E_{\alpha }(-t)=\int_{0}^{\infty }e^{-t\tau }\,d\nu _{\alpha }(\tau ),\quad
\forall \,t\geq 0  \label{F1}
\end{equation}%
for $\nu _{\alpha }$ being the probability measure on $\mathbb{R}_{0}^{+}$ 
\begin{equation}
d\nu _{\alpha }\left( \tau \right) :=\alpha ^{-1}\tau ^{-1-1/\alpha
}f_{\alpha }(\tau ^{-1/\alpha })\,d\tau  \label{A2}
\end{equation}%
where $f_{\alpha }$ is the $\alpha $-stable probability density given by 
\begin{equation*}
\int_{0}^{\infty }e^{-t\tau }f_{\alpha }\left( \tau \right) \,d\tau
=e^{-t^{\alpha }},\quad 0<\alpha <1
\end{equation*}

The complete monotonicity property and the integral representation (\ref{F1}%
) of $E_{\alpha }$ may be extended to complex arguments.

\begin{lemma}
\label{Aug09-eq3}For any $z\in \mathbb{C}$ such that $\mathrm{Re}(z)\geq 0$,
the following representation holds 
\begin{equation*}
E_{\alpha }(-z)=\int_{0}^{\infty }e^{-z\tau }\,d\nu _{\alpha }(\tau ),\quad
0<\alpha \leq 1
\end{equation*}
\end{lemma}

\noindent \textbf{Proof.} According to \cite{Pollard}, for each $0<\alpha<1$
fixed, for all $t\geq 0$ one has 
\begin{eqnarray}
E_{\alpha }(-t)&=&\int_{0}^{\infty }e^{-t\tau }\,d\nu _{\alpha }(\tau ), 
\notag \\
&=& \sum_{n=0}^{\infty }\frac{(-t)^{n}}{n!}\int_{0}^{\infty }\tau ^{n}\,d\nu
_{\alpha }(\tau )  \label{Sep09-eq5}
\end{eqnarray}
Comparing (\ref{Sep09-eq5}) with the Taylor expansion (\ref{1.9}) of $%
E_\alpha$, one concludes that the moments of the measure $\nu_\alpha$ are
given by 
\begin{equation*}
m_{n}(\nu _{\alpha }):=\int_{0}^{\infty }\tau ^{n}\,d\nu _{\alpha }(\tau )=%
\frac{n!}{\Gamma (\alpha n+1)},\quad n\in \mathbb{N}_{0}
\end{equation*}

For complex values $z$ let 
\begin{equation*}
I(-z):=\int_{0}^{\infty }e^{-z\tau }\,d\nu _{\alpha }(\tau )
\end{equation*}%
which is finite provided $\mathrm{Re}(z)\geq 0$. For each $z\in \mathbb{C}$
such that $\mathrm{Re}(z)\geq 0$ one then obtains 
\begin{equation*}
I(-z)=\sum_{n=0}^{\infty }\frac{\left( -z\right) ^{n}}{n!}\left(
\int_{0}^{\infty }\tau ^{n}\,d\nu _{\alpha }(\tau )\right)
=\sum_{n=0}^{\infty }\frac{\left( -z\right) ^{n}}{n!}m_{n}(\nu _{\alpha
})=\sum_{n=0}^{\infty }\frac{\left( -z\right) ^{n}}{\Gamma (\alpha n+1)}%
=E_{\alpha }(-z)
\end{equation*}%
leading to the integral representation 
\begin{equation*}
E_{\alpha }(-z)=\int_{0}^{\infty }e^{-z\tau }d\nu _{\alpha }(\tau )
\end{equation*}%
for all $z\in \mathbb{C}$ such that $\mathrm{Re}(z)\geq 0$.\hfill $%
\blacksquare \medskip $


\begin{thebibliography}{99}
\bibitem{Albeverio1} S. Albeverio, Yu. G. Kondratiev and M. R\"{o}ckner; 
\textit{Analysis and geometry on configuration spaces}, J. Funct. Anal. 154
(1998) 444--500.

\bibitem{Albeverio2} S. Albeverio, Yu. G. Kondratiev and M. R\"{o}ckner; 
\textit{Analysis and geometry on configuration spaces: The Gibbsian case},
J. Funct. Anal. 157 (1998) 242--291.

\bibitem{FCAA} F. Cipriano, H. Ouerdiane and R. Vilela Mendes; \textit{%
Stochastic solution of a KPP-type nonlinear fractional differential equation}%
, Fract. Calcul. and Appl. Analysis 12 (2009) 47--57.

\bibitem{Fel71} W. Feller; \textit{An introduction to probability theory and
its applications, vol. II}, Second Edition, John Wiley \& Sons, N. Y. 1971.

\bibitem{FicWin93} K.-H. Fichtner and G. Winkler; \textit{Generalized
Brownian motion, point processes and stochastic calculus for random fields},
Math. Nachr. 161 (1993) 291--307.

\bibitem{GV} I. M. Gelfand and Ya. N. Vilenkin; \textit{Generalized
Functions, vol. IV}, Academic Press, N. Y. and London 1968.

\bibitem{YuKu} Yu. G. Kondratiev and T. Kuna; \textit{Harmonic analysis on
configuration space I. General theory}, Infin. Dimens. Anal. Quantum Probab.
Relat. Top. 5 (2002) 201--233.

\bibitem{KoKu100} Yu. G. Kondratiev, T. Kuna and M. J. Oliveira; \textit{%
Analytic aspects of Poissonian white noise analysis}, Methods Funct. Anal.
Topology 8 (2002) 15--48.

\bibitem{KKO04} Yu. G. Kondratiev, T. Kuna and M. J. Oliveira \textit{On the
relations between Poissonian white noise analysis and harmonic analysis on
configuration spaces}, J. Funct. Anal. 213 (2004) 1--30.

\bibitem{Le72} A. Lenard; \textit{Correlation functions and the uniqueness
of the state in classical statistical mechanics}, Commun. Math. Phys. 30
(1973) 35--44.

\bibitem{Le75a} A. Lenard; \textit{States of classical statistical
mechanical systems of infinitely many particles I}, Arch. Rational Mech.
Anal. 59 (1975) 219--239.

\bibitem{Le75b} A. Lenard; \textit{States of classical statistical
mechanical systems of infinitely many particles II}, Arch. Rational Mech.
Anal. 59 (1975) 241--256.

\bibitem{Mainardi1} F. Mainardi, R. Gorenflo and E. Scalas; \textit{A
fractional generalization of the Poisson process}, Vietnam J. of Mathematics
32 (2004) 53--64.

\bibitem{Oliveira1} M. J. Oliveira; \textit{Configuration Space Analysis and
Poissonian White Noise Analysis}, Ph. D. Thesis, Lisbon 2002,
http://www.math.uni-bielefeld.de/igk/study-materials/corpo\_2.pdf.

\bibitem{P67} K. R. Parthasarathy; \textit{Probability Measures on Metric
Spaces}, Academic Press, N. Y. 1967.

\bibitem{Pollard} H. Pollard; \textit{The complete monotonic character of
the Mittag-Leffler function }$E_{\alpha }\left( -x\right) $, Bull. Amer.
Math. Soc. 54 (1948) 1115--1116.

\bibitem{S94} H. Shimomura; \textit{Poisson measures on the configuration
space and unitary representations of the group of diffeomorphisms}, J. Math.
Kyoto Univ. 34 (1994) 599--614.
\end{thebibliography}
\end{document}